\documentclass[a4paper,12pt]{article}
\usepackage{graphicx}
\usepackage{a4wide}
\usepackage[english]{babel}
\usepackage{amsfonts}
\usepackage{amsmath}
\usepackage{amssymb}
\usepackage{dsfont}
\newtheorem{lemma}{Lemma}

\newtheorem{theorem}{Theorem}

\allowdisplaybreaks

\date{}

\newcommand {\E} {\mathbb{E}}
\DeclareMathOperator {\var}{Var_{[0,1]}}

\newcommand {\p} {\mathbb{P}}
\newcommand {\Z} {\mathbb{Z}}
\newcommand {\N} {\mathbb{N}}
\newcommand {\R} {\mathbb{R}}

\newcommand {\ve} {\varepsilon}

\makeatletter\def\blfootnote{\xdef\@thefnmark{}\@footnotetext}\makeatother

\begin{document}
\title{\bf On the law of the iterated logarithm for permuted lacunary sequences}
\author{C.\ Aistleitner\footnote{Graz University of Technology, Department for Analysis and Computational Number Theory, Steyrergasse 30, 8010 Graz, Austria \mbox{e-mail}: \texttt{aistleitner@math.tugraz.at}.
Research supported by FWF grant S9603-N23.}, I.\ Berkes\footnote{ Graz University of Technology, Institute of
Statistics, M\"unzgrabenstra{\ss}e 11, 8010 Graz, Austria.  \mbox{e-mail}: \texttt{berkes@tugraz.at}. Research
supported by the FWF Doctoral Program on Discrete Mathematics (FWF DK W1230-N13), FWF grant S9603-N23 and OTKA grants K 67961 and K 81928.} and R.\ Tichy\footnote{Graz University of Technology, Department for Analysis and Computational Number Theory, Steyrergasse 30, 8010 Graz, Austria. \mbox{e-mail}: \texttt{tichy@tugraz.at}.
Research supported by the FWF Doctoral Program on Discrete Mathematics (FWF DK W1230-N13) and FWF grant S9603-N23.}}

\maketitle

\emph{Dedicated to the memory of Professor Anatolii Alexeevitch Karatsuba.}

\abstract{It is known that for any smooth periodic function $f$ the sequence $(f(2^kx))_{k\ge 1}$ behaves like a sequence of i.i.d.\ random variables, for example, it satisfies the central limit theorem and the law of the iterated logarithm. Recently Fukuyama showed that permuting $(f(2^kx))_{k\ge 1}$ can ruin the validity of the law of the iterated logarithm, a very surprising result. In this paper we present an optimal condition on $(n_k)_{k\ge 1}$, formulated in terms of the number of solutions of certain Diophantine equations, which ensures the validity of the law of the iterated logarithm for any permutation of the sequence $(f(n_k x))_{k \geq 1}$. A similar result is proved for the discrepancy of the sequence $(\{n_k x\})_{k \geq 1}$, where $\{ \cdot \}$ denotes fractional part.}


\section{Introduction}

Given a sequence  $(x_1, \dots, x_N)$ of real numbers, the value
$$
D_N = D_N (x_1, \dots, x_N) = \sup_{0 \leq a < b < 1} \left|
\frac{\sum_{k=1}^N \mathds{1}_{[a,b)} (\{x_k\})}{N} - (b-a)
\right|
$$
is called the discrepancy of the sequence. Here
$\mathds{1}_{[a,b)}$  denotes the indicator function of the
interval $[a,b)$ and $\{ \cdot \}$ denotes fractional part. An
infinite sequence $(x_n)_{n\ge 1}$ is called uniformly distributed
mod 1 if $D_N(x_1, \ldots x_N) \to 0$ as $N\to \infty$. Weyl
\cite{we} proved that for any increasing sequence $(n_k)_{k\ge 1}$ of
integers, $(n_kx)_{k\ge 1}$ is uniformly distributed mod 1 for almost all $x\in{\mathbb R}$ in the sense of the Lebesgue
measure. Computing the discrepancy of this sequence is a
difficult problem and the precise asymptotics is known only in a
few cases. Philipp \cite{plt} proved that if $(n_k)_{k\ge 1}$ satisfies the
Hadamard gap condition
\begin{equation} \label{had}
n_{k+1}/n_k \geq q > 1 \qquad (k=1, 2, \ldots),
\end{equation}
then the discrepancy of $(\{n_kx\})_{k \geq 1}$ obeys the law of the iterated
logarithm, i.e.
\begin{equation}\label {phlil}
\frac{1}{4 \sqrt{2}} \leq \limsup_{N \to \infty} \frac{N D_N(n_k
x)} {\sqrt{2 N \log \log N}} \leq D(q) \quad \textup{a.e.},
\end{equation}
where $D(q)$ is a number depending on $q$. (For the simplicity of the
notation, we will write $D_N(x_k)$ instead of $D_N(x_1, \ldots
x_N)$.) Note that for the discrepancy of an i.i.d.\ nondegenerate sequence $(X_k)_{k\ge 1}$ we have
\begin{equation}\label{chs}
\limsup_{N \to \infty} \frac{N D_N(X_k)} {\sqrt{2 N \log \log N}}
=1/2 \quad \textup{a.s.}
\end{equation}
by the Chung-Smirnov law of the iterated logarithm (see e.g. \cite[p.\ 504]{sw}). A comparison of (\ref{phlil}) and (\ref{chs}) shows that for Hadamard lacunary $(n_k)_{k\ge 1}$, the
sequence $(\{n_kx\})_{k \geq 1}$ of functions on $(0, 1)$  behaves like a
sequence of i.i.d.\ random variables. (Note, however, that for certain values of $x$ the distribution of $(\{n_k x\})_{k \geq 1}$ can differ significantly from the uniform distribution even under (\ref{had}), see e.g. \cite{peres}.) Interestingly, the analogy between lacunary sequences and sequences of i.i.d. random variables is not complete. Fukuyama \cite{ft} determined the limsup in (\ref{phlil}) for the sequences $n_k=\theta^k$, $\theta>1$; his results show that the limsup is different from $1/2$ for any integer $\theta\ge 2$. (On the other hand, in \cite{ft} it is shown that if $\theta^r$ is irrational for $r=1, 2, \ldots$, then the limsup in (\ref{phlil}) equals to the i.i.d.\ value $1/2$.)  Aistleitner \cite{ai2} constructed a Hadamard lacunary sequence $(n_k)_{k\ge 1}$ such that the limsup in (\ref{phlil}) is not a constant a.e.\ and Fukuyama and Miyamoto \cite{fuku2^n-1} showed that this actually happens for $n_k=2^k-1$. Even more surprisingly, Fukuyama \cite{ft2} showed that even if the limsup in (\ref{phlil}) is a constant (e.g., for $n_k=2^k)$, the value of the limsup can change by permuting the sequence $(n_k)_{k\ge 1}$, a phenomenon radically different from i.i.d.\ behavior, which is clearly permutation-invariant.\\

The previous results show that the behavior of $(\{n_kx\})_{k\ge 1}$ is very delicate, exhibiting both probabilistic and number-theoretic phenomena. It is natural to ask for which $(n_k)_{k\ge 1}$ the behavior of  $(\{n_kx\})_{k\ge 1}$ follows exactly i.i.d.\ behavior, for example, when is the limsup in (\ref{phlil}) equal to $1/2$ a.e and under what conditions is the value of the limsup permutation-invariant. A near optimal number-theoretic condition for $\limsup=1/2$ a.e. was given by Aistleitner \cite{ao} and the purpose of the present paper is to give an optimal condition for the permutation-invariance of the LIL.\\

Let $f$ be a measurable function satisfying
\begin{equation}\label{fcond1}
f(x+1)=f(x), \quad \int_0^1 f(x)dx=0, \quad \text{Var}_{[0, 1]} (f)<\infty.
\end{equation}
A profound study of the behavior of the sequence $(f(n_kx))_{k\ge 1}$ for Hadamard lacunary $(n_k)_{k\ge 1}$ was given in Gaposhkin \cite{gap1966, gap1970}. By a classical theorem of Kac \cite{kac47}, under (\ref{fcond1})  the sequence $(f(2^kx))_{k\ge 1}$ satisfies the CLT and Erd\H{o}s  and Fortet showed (see \cite[p.\ 646]{kac49}) that the CLT generally fails for $(f((2^k-1)x))_{k\ge 1}$ (see also \cite{conze}). Gaposhkin showed that $(f(n_kx))_{k\ge 1}$ satisfies the central limit theorem provided  $n_{k+1}/n_k$ is an integer for all $k\ge 1$ or if $n_{k+1}/n_k\to\alpha$, where $\alpha^r$ is irrational for $r=1, 2, \ldots$. More generally, he showed that $(f(n_kx))_{k\ge 1}$ satisfies the central limit theorem
provided that for any nonzero integers $a, b, c$ the number of solutions of the Diophantine equation
\begin{equation} \label{gap}
a n_k + b n_l = c, \qquad 1 \le k,l \leq N
\end{equation}
is bounded by a constant $K(a, b)$, independent of $c$.
Aistleitner and Berkes \cite{aibe} showed that the CLT remains valid if for any nonzero integers $a, b, c$ the number of solutions of (\ref{gap}) is $o(N)$, uniformly in $c$, and this condition is best possible. Aistleitner \cite{ao} also proved that replacing $o(N)$ by $O(N/(\log N)^{1+\varepsilon})$ in the previous theorem, the limsup in (\ref{phlil}) equals $1/2$. As we will see, a two-term Diophantine condition will also give the precise condition for the permutation-invariance of the LIL for $D_N(n_kx)$.

\section{Results}

In what follows, we write $\|f\|$ for the $L^2(0,1)$ norm of a function $f$.

\begin{theorem} \label{thlils}
Let $(n_k)_{k \geq 1}$ be a sequence of positive integers satisfying (\ref{had}), such that for any fixed integers $a\ne 0$, $b\ne 0$, $c$ the number of
solutions of the Diophantine equation (\ref{gap})
is bounded by a constant $K(a, b)$ independent of $c$, where for $c=0$ we require
also $k\ne l$. Let $f$ be a function satisfying (\ref{fcond1}).
Then for any permutation $\sigma: \N \to \N$ we have 
$$
\limsup_{N\to \infty} \frac{\sum_{k=1}^N
f(n_{\sigma (k)} x)}{\sqrt{2N\log\log N}} = \|f\| \quad\textup{a.e.}
$$
\end{theorem}

As a consequence of Theorem \ref{thlils} we obtain the following metric discrepancy result.

\begin{theorem}\label{dlil}
Let $(n_k)_{k\ge 1}$ satisfy the assumptions of Theorem 1. Then for any permutation $\sigma: {\mathbb N}\to {\mathbb N}$ we have 
\begin{equation}\label {lilinfgap2} \limsup_{N \to \infty} \frac{ND_N(n_{\sigma(k)}x)}{\sqrt{2 N \log \log N}} =\frac{1}{2} \quad
\textup{a.e.}
\end{equation}
\end{theorem}

At the end of our paper we will show that if there exist  integers $a\ne 0, b\ne 0$ and $c$ such that the Diophantine equation (\ref{gap}) has infinitely many solutions $(k,l), k \neq l$, then the conclusion of Theorem \ref{thlils} fails to hold for appropriate $f$. In fact, we can even obtain a non-constant limsup, which perfectly matches the results in \cite{ai, ai2}. This shows that the Diophantine condition in Theorem \ref{thlils} is essentially optimal.\\

We stress that in Theorems \ref{thlils} and \ref{dlil} we bounded the number of solutions of (\ref{gap}) also for $c=0$ and thus $n_k=2^k$ does not satisfy this condition. In fact, the conclusion of both theorems is false for $n_k=2^k$:  with the identity permutation $\sigma(k)=k$ the limsup in Theorem \ref{thlils} equals
$$
\left(\|f\|^2+2\sum_{k=1}^\infty \int_0^1 f(x)f(2^kx)\, dx\right)^{1/2}
$$
(see \cite{iz,ma}) and the limsup in Theorem \ref{dlil}
is $\sqrt{42}/9$ by the theorem of Fukuyama \cite{ft}. As we mentioned in the Introduction, for the CLT with a limit distribution of unspecified variance, it suffices to bound the number of solutions of (\ref{gap}) for coefficients $a, b, c$ all different from 0.\\

As the proof of Theorem \ref{thlils} will show, in the case when $f$ is a trigonometric polynomial
of degree $d$, it suffices to assume the Diophantine condition only with coefficients $a,b$ satisfying  $1 \leq |a|\leq d, 1 \leq |b| \leq d$. In particular, in the trigonometric case $f(x)=\cos 2\pi x$ it suffices  to allow only coefficients $\pm 1$, when the Diophantine condition in Theorem \ref{thlils} is satisfied for any Hadamard lacunary sequence $(n_k)_{k\ge 1}$
(see e.g.\ Zygmund \cite[pp.~203-204]{zt}). Thus we obtain the following corollary of Theorem \ref{thlils}, which is a permutation invariant version of the Erd\H os-G\'al LIL in \cite{eg}.

\begin{theorem}\label{permcltlil2}
Let $(n_k)_{k \geq 1}$ be a sequence of positive integers satisfying
the Hadamard gap condition (\ref{had}), and let $\sigma: \N \to \N$ be a permutation of the set of positive integers. Then
\begin{equation}\label{permlil}
\limsup_{N\to\infty} \frac{\sum_{k=1}^N \cos 2\pi
n_{\sigma(k)}x}{\sqrt{2 N \log \log N}}= \frac{1}{\sqrt{2}} \qquad \textup{a.e.}
\end{equation}
\end{theorem}

If we assume the slightly stronger gap condition
\begin{equation}\label{infgap}
n_{k+1}/n_k\to\infty
\end{equation}
instead of Hadamard's gap condition (\ref{had}), then the behavior of $f(n_kx)$ is permutation-invariant, regardless
the number theoretic structure of $(n_k)_{k \geq 1}$. In fact, any sequence $(n_k)_{k \geq 1}$ satisfying the
gap condition (\ref{infgap}) satisfies the Diophantine condition in Theorem \ref{thlils} automatically. This follows  from the fact
that for arbitrary fixed nonzero integers $a, b$ the set-theoretic union of the sequences $(a n_k)_{k \ge 1}$ and $(b n_k)_{k\ge 1}$,
arranged in increasing order, satisfies the Hadamard gap condition (\ref{had}) and consequently the Diophantine condition in Theorem \ref{thlils}. Thus Theorem \ref{thlils} implies the following

\begin{theorem}\label{th1}
Let $(n_k)_{k \geq 1}$ be a sequence of positive integers satisfying the gap
condition (\ref{infgap}). Then for any permutation $\sigma:
{\mathbb N}\to {\mathbb N}$ of the integers and for any function $f$ satisfying (\ref{fcond1}) we have
\begin{equation}\label{lilperm1}
\limsup_{N\to \infty} \frac{\sum_{k=1}^N
f(n_{\sigma (k)} x)}{\sqrt{2N\log\log N}} =\|f\| \quad\textup{a.e.}\\
\end{equation}
Moreover, for any permutation $\sigma$ of $\mathbb{N}$ we have
\begin{equation}\label {lilinfgap} \limsup_{N \to \infty} \frac{ND_N(n_{\sigma(k)}x)}{\sqrt{2 N \log \log N}} =\frac{1}{2} \quad
\textup{a.e.}
\end{equation}
\end{theorem}

\section{Proofs}


To deduce Theorem \ref{dlil} from Theorem \ref{thlils}, we write
$\mathbf{I}_{[a,b)}$ for the indicator of the interval $[a,b)$,
centered at expectation and extended with period 1, i.e.
$$
\mathbf{I}_{[a,b)}(x) = \mathds{1}_{[a,b)} (\langle x \rangle) -
(b-a),
$$
where $\langle \cdot \rangle$ denotes the fractional part. For lacunary sequences (and also for
permutations of lacunary sequences)
$$
\limsup_{N \to \infty} \sup_{0 \leq a < b \leq 1} \frac{ \left|
\sum_{k=1}^N \mathbf{I}_{[a,b)} (n_k x) \right|}{\sqrt{2 N \log
\log N}} = \sup_{0 \leq a < b \leq 1} \limsup_{N \to \infty}
\frac{ \left| \sum_{k=1}^N \mathbf{I}_{[a,b)} (n_k x)
\right|}{\sqrt{2 N \log \log N}}\quad \textup{a.e.}
$$
(see \cite[Theorem 1]{fukusmall}). Thus Theorem \ref{thlils} yields
\begin{eqnarray*}
& & \limsup_{N \to \infty} \frac{N D_N(n_{\sigma(k)} x)}{\sqrt{2 N \log
\log N}} \\
& = & \sup_{0 \leq a < b \leq 1} \limsup_{N \to \infty}
\frac{ \left| \sum_{k=1}^N \mathbf{I}_{[a,b)} (n_{\sigma(k)} x)
\right|}{\sqrt{2 N \log \log N}} = \sup_{0 \leq a < b \leq 1}
\left\|\mathbf{I}_{[a,b)} \right\| = \frac{1}{2}
\end{eqnarray*}
almost everywhere.\\

We now turn to the proof of Theorem \ref{thlils}.
In the sequel we will assume that the function $f(x)$, the sequence $(n_k)_{k \geq 1}$
and the permutation $\sigma$ are fixed. We assume that sequence $(n_k)_{k \geq 1}$ is a lacunary sequence satisfying the Diophantine condition in Theorem \ref{thlils}. The method which we use for the proof of Theorem \ref{thlils} is a multidimensional version of the method which was used in to prove the permutation-invariant CLT in our paper \cite{abt}. We recommend \cite{abt} as an introduction to the methods which are used in the present paper.\\

We begin with some auxiliary results.
\begin{lemma} \label{rp}
Let $A_1, A_2, \dots$ be arbitrary events, satisfying
$$
\sum_{m=1}^\infty \p (A_m) = \infty
$$
and
$$
\liminf_{N \to \infty} \frac{\sum_{n=1}^N \sum_{m=1}^N \p (A_m
A_n)}{\left( \sum_{m=1}^N \p (A_m) \right)^2} =1.
$$
Then
$$
\p \left( \prod_{n=1}^\infty \sum_{m=n}^\infty A_m \right) = 1,
$$
i.e. with probability 1 infinitely many of the events
$A_m$ occur.
\end{lemma}
A proof of this lemma can be found in \cite{rp}.\\

\begin{lemma} \label{lemmacf}
Let $P_1,P_2$ be probability measures on $\R^2$, and write $p_1,
p_2$ for the corresponding characteristic functions. Then for all $T_1,
T_2,\delta_1,\delta_2,x,y>0$
\begin{eqnarray*}
& & \left| P_1^*([-x,x] \times [-y,y]) - P_2^*([-x,x] \times [-y,y]) \right| \\
& \leq & xy \int_{(s,t) \in
[-T_1,T_1] \times [-T_2,T_2]} |p_1(s,t)-p_2(s,t)| ~d(s,t) \\
& & + xy  \left(  \delta_1^{-1} \delta_2~ \exp\left(- T_1^2
\delta_1^2/2 \right) + \delta_1 \delta_2^{-1} ~\exp \left(-
T_2^2 \delta_2^2/2 \right) \right),
\end{eqnarray*}
where
$$
P_1^* = P_1 \star H, \quad P_2^* = P_2 \star H,
$$
and $H$ is a two-dimensional normal distribution with density
$$
(2 \pi \delta_1 \delta_2)^{-1} ~e^{-\delta_1^2 u^2/2 -
\delta_2^2 v^2/2}.
$$
\end{lemma}
\emph{Proof:~} Assume that $T_1, T_2,\delta_1,\delta_2
>0$ are fixed. Letting
$$
h(s,t) = e^{-\delta_1^2 s^2/2 - \delta_2^2 t^2 /2}
$$
denote the characteristic function of $H$, we have $p_1^*=p_2 h$
and $p_2^*=p_2 h$. Writing $\gamma_1$ and $\gamma_2$ for the
densities of $P_1^*$ and $P_2^*$, respectively, we have
\begin{eqnarray*}
|\gamma_1(u,v) - \gamma_2(u,v)| & \leq & (2 \pi)^{-2} \left|
\int_{\R^2} e^{-isu-itv} \left( p_1^*(s,t)-p_2^*(s,t) \right)
~d(s,t)
\right| \\
& \leq & (2 \pi)^{-2} \int_{\R^2} \left| p_1(s,t)-p_2(s,t)\right|
|h(s,t)|~d(s,t) \\
& \leq & (2 \pi)^{-2} \int_{(s,t) \in [-T_1,T_1] \times
[-T_2,T_2]} |p_1(s,t)-p_2(s,t)| ~d(s,t) \\
& &  + (2 \pi)^{-2} ~2 \int_{(s,t) \not\in [-T_1,T_1] \times
[-T_2,T_2]} |h(s,t)| ~d(s,t).
\end{eqnarray*}
Therefore
\begin{eqnarray*}
& & \left| P_1^* ([-x,x] \times[-y,y]) - P_2^*([-x,x] \times [-y,y]) \right| \\
& \leq & \int_{[-x,x] \times [-y,y]} |\gamma_1(u,v) -
\gamma_2(u,v)| ~d(u,v) \\
& \leq & 4xy (2 \pi)^{-2} \int_{(s,t) \in [-T_1,T_1] \times
[-T_2,T_2]} |p_1(s,t)-p_2(s,t)| ~d(s,t) \\
& & + 4 xy (2 \pi)^{-2} ~2 \int_{(s,t) \not\in [-T_1,T_1] \times
[-T_2,T_2]} |h(s,t)| ~d(s,t).
\end{eqnarray*}
Now
\begin{eqnarray*}
& & \int_{(s,t) \not\in [-T_1,T_1] \times [-T_2,T_2]} |h(s,t)| ~d(s,t) \\
& = & \int_{(s,t) \not\in [-T_1,T_1] \times [-T_2,T_2]}
e^{-\delta_1^2 s^2/2 - \delta_2^2 t^2 /2} ~d(s,t) \\
& \leq & (2\pi)^{-1}\delta_1 \delta_2 - \left(
(2\pi)^{-1/2}\delta_1 - 4 \delta_1^{-1} \exp \left(-
T_1^2 \delta_1^2/2 \right) \right) \times \nonumber\\
& & \quad \times  \left(
(2\pi)^{-1/2}\delta_2 - 4\delta_2^{-1} \exp
\left( - T_2^2 \delta_2^2/2 \right) \right)  \\
& \leq & \delta_1^{-1} \delta_2~ \exp\left(- T_1^2
\delta_1^2/2 \right) + \delta_1 \delta_2^{-1} ~\exp \left(-
T_2^2 \delta_2^2/2 \right)
\end{eqnarray*}
proves the lemma.\qquad $\square$\\

The following lemma is a one-dimensional version of Lemma \ref{lemmacf}, which can be shown in the same way (a proof is contained in \cite{abt}).
\begin{lemma}[{\cite[Lemma 4.2]{abt}}] \label{lemmaclt2} \label{lemmacf1dim}
Let $P_1,P_2$ be probability measures on $\R$, and write $p_1,
p_2$ for the corresponding characteristic functions. Let
$$
P_1^* = P_1 \star H,\quad P_2^* = P_2 \star H,
$$
where $H$ is a normal distribution with mean zero and standard deviation $\delta$. Then for all $T>0$
\begin{eqnarray*}
| P_1^*([-x,x]) - P_2^*([-x,x]) | & \leq & x \int_{s \in
[-T,T]} |p_1(s)-p_2(s)| ~ds + 4 x \delta^{-1} e^{- T^2 \delta^2/2}.
\end{eqnarray*}
\end{lemma}

\begin{lemma} \label{lemmacf1dimstern}
Let $P_1,P_2$ be probability measures on $\R$, and write $p_1,
p_2$ for the corresponding characteristic functions. Then for any $S,T,\delta>0$
\begin{eqnarray*}
& & \left| P_1 ([-x-S,x+S]) - P_2 ([-x+S,x-S]) \right|\\
& \leq & x \int_{s \in  [-T,T]}
|p_1(s)-p_2(s)| ~ds + 4 x \delta^{-1} e^{- T^2 \delta^2/2} + 4 e^{- S^2/(2 \delta^2)}.
\end{eqnarray*}
\end{lemma}
\emph{Proof:~} Let $S,T,\delta$ be fixed. Let $H$ be a normal
distribution with mean zero and variance $\delta^2$, and set $P_1^*= P_1 \star H$ and $P_2^* = P_2 \star H$. Then
\begin{eqnarray*}
& & \left| P_1 ([-x-S,x+S]) - P_2 ([-x+S,x-S]) \right| \\
& \leq & \left| P_1^* ([-x,x]) - P_2^*([-x,x]) \right| + 2 H
\left( \R \backslash [-S,S]\right).
\end{eqnarray*}
Now
\begin{eqnarray*}
2 H(\R \backslash [-S,S]) & \leq & \frac{2}{\delta \sqrt{2 \pi}} \int_{s \not\in [-S,S]} e^{-s^2/(2 \delta^2)} ~ds \\
& \leq & 4 e^{- S^2/(2 \delta^2)},
\end{eqnarray*}
and thus by Lemma \ref{lemmaclt2}
\begin{eqnarray*}
& & \left| P_1 ([-x-S,x+S]) - P_2 ([-x+S,x-S]) \right|\\
& \leq & x \int_{s \in
[-T,T]} w(s) ~ds + 4 e^{- S^2/(2 \delta^2)} + 4 x \delta^{-1} e^{- T^2 \delta^2/2},
\end{eqnarray*}
which proves Lemma \ref{lemmacf1dimstern}. \qquad $\square$\\

Let now $\theta>1$, $\ve>0$ and $d\ge 1$ be arbitrary, but fixed. To simplify notations we assume in the sequel that $f$ is even. We will also assume that
$\|f\|>0$, since otherwise Theorem \ref{thlils} is trivial. We write
$$
f(x) \sim \sum_{j=1}^\infty a_j \cos 2 \pi j x = p(x)+r(x),
$$
where
$$
p(x)=\sum_{j=1}^d a_j \cos 2 \pi j x, \quad
r(x)=\sum_{j=d+1}^\infty a_j \cos 2 \pi j n_k x.
$$

\begin{lemma}[{\cite[Lemma 3.1]{ao}}] \label{lemmar}
$$
\limsup_{N \to \infty} \frac{\sum_{k=1}^N r(n_{\sigma(k)}
x)}{\sqrt{2 N \log \log N}} \leq C d^{-1/4} \quad \textup{a.e.}
$$
for some number $C$ which is independent of $d$ and $\sigma$.
\end{lemma}

This lemma is valid for general lacunary sequences without any condition on the number of
solutions of Diophantine equations. It is easy to see that the
proof of \cite[Lemma 3.1]{ao} is valid not only for strictly
increasing lacunary sequences, but also for permutations of
lacunary sequences, since it is based on a result of Philipp
\cite{plt}, which also has this property. We also recall that
$\var f \leq 1$ implies $|a_j| \leq j^{-1}, ~j \geq 1$
(cf. Zygmund \cite[p.~48]{zt}), and
$$
\|r\|^2 \leq \sum_{j=d+1}^\infty |a_j|^2 \leq \sum_{j=d+1}^\infty
j^{-2} \leq d^{-1}.
$$

\begin{lemma} \label{lemmap}
$$
\limsup_{N \to \infty} \frac{\sum_{k=1}^N p(n_{\sigma(k)}
x)}{\sqrt{2 N \log \log N}} = \|p\| \quad \textup{a.e.}
$$
\end{lemma}

In case $\|p\|=0$ Lemma \ref{lemmap} is trivial. Therefore, to
simplify formulas, we will assume in the sequel, without loss of generality, that $\|p\|=1$, and prove
$$
\limsup_{N \to \infty} \frac{\sum_{k=1}^N p(n_{\sigma(k)}
x)}{\sqrt{2 N \log \log N}} = 1 \quad \textup{a.e.}
$$
Since a finite number of elements of $(n_k)_{k \geq 1}$ does not influence the
asymptotic behavior of $(n_{\sigma(k)} x)_{k \geq 1}$ we can also assume that
\begin{equation} \label{j1j2}
a n_k + b n_l = 0
\end{equation}
does not have any solution $(k,l),~k \neq l$ for $1 \leq |a| \leq d,1 \leq |b|
\leq d$.\\

Since $d$ can be chosen arbitrarily large, Lemma \ref{lemmar} and and Lemma \ref{lemmap} together imply Theorem \ref{thlils}.
Therefore it remains to prove Lemma \ref{lemmap}. The proof of
this lemma is crucial, and will be given in two parts below. The
main ingredient is Lemma \ref{2dim}, which is formulated and
proven below. For the proof we use ideas of R{\'e}v{\'e}sz \cite{rt}.\\

We define
\begin{eqnarray}
\mu_k & = & n_{\sigma(k)}, \qquad k \geq 1 \nonumber\\
\Delta_m^* & = & \left\{ k \geq 1: ~ \theta^m
\leq k < \theta^{m+1} \right\}, \qquad m \geq 1. \nonumber
\end{eqnarray}
We rearrange the sequence $(\mu_k)_{k \geq 1}$ in such a way that it is
increasing within the blocks $\Delta_m^*$ and
call this new sequence $(\nu_k)_{k \geq 1}$. In other words, $(\nu_k)_{k \geq
1}$ consists of the same elements as $(\mu_k)_{k \geq 1}$ (and
$(n_{\sigma(k)})_{k \geq 1}$), but satisfies
$$
\nu_k < \nu_l \quad \textrm{if} \quad k<l \quad \textrm{and} \quad
k,l \in \Delta_m^* \quad \textrm{for some}
\quad m \geq 1.
$$
Moreover, we define
\begin{eqnarray}
\overline{\Delta}_m & = & \left\{ k \in
\Delta_m^*:~ \nexists l \in \bigcup_{n=1}^{m -
\log_\theta m} \Delta_n^*:~
\frac{\nu_k}{\nu_l} \in
\left[ \frac{1}{2d},2d \right] \right\},\quad m \geq 1, \label{constro}\\
\Delta_m & = & \Big\{k \in \overline{\Delta}_m:~ \left(k \mod
\left( \left\lceil \sqrt{m} + \log_q (2 d) \right\rceil \right)
\right) \nonumber\\
& & \qquad \not\in \left\{0,\dots,\left\lceil \log_q (2d)\right\rceil
\right\} \Big\}, \quad m \geq 1,
\nonumber\\
\Delta_m^{(h)} & = & \left\{ k \in \Delta_m:~ \frac{k}{\lceil
\sqrt{m}+\log_q (2 d) \rceil}
\in [h,h+1) \right\}, \quad h \geq 0, m \geq 1, \nonumber\\
\eta_m & = & \frac{\sum_{k \in \Delta_m} p (\nu_k
x)}{\sqrt{|\Delta_m|}}, \quad m \geq 1, \nonumber\\
\alpha_m (s) & = & \prod_{h \geq 0} \left(1+\frac{i s \sum_{k \in
\Delta_m^{(h)}} \sum_{j=1}^d a_j \cos(2 \pi j \nu_k
x)}{\sqrt{|\Delta_m|}} \right),
\quad m \geq 1, \nonumber\\
\beta_m & = & \sum_{h \geq 0} ~\sum_{k_1,k_2 \in \Delta_m^{(h)}}
\sum_{j_1,j_2=1}^d \frac{a_{j_1} a_{j_2}}{2} \cos (2 \pi (j_1
\nu_{k_1} + j_2
\nu_{k_2}) x ) \nonumber\\
& & + \sum_{h \geq 0}~\underbrace{\sum_{k_1,k_2 \in
\Delta_m^{(h)}} \sum_{j_1,j_2=1}^d}_{(k_1,j_1) \neq (k_2,j_2)}
\frac{a_{j_1} a_{j_2}}{2} \cos (2 \pi
(j_1 \nu_{k_1} - j_2 \nu_{k_2}) x ), \quad m \geq 1,\nonumber\\
\varphi_{m,n} (s,t)& = & \E \left( e^{i s \eta_m +
i t \eta_n} \right), \quad m,n \geq 1, \quad s,t \in \R \nonumber
.
\end{eqnarray}
Here $|\Delta_m|$ denotes the number of elements of the set
$\Delta_m$. Throughout the paper $\log x$ will be understood as
$\max\{1,\log x\}$. $\sum_{k=m}^n$ means $\sum_{m \leq k \leq n}$,
if $m,n$ are not integers.\\

We observe that for $m \geq 1$
\begin{eqnarray*}
|\overline{\Delta}_m| & \geq & |\Delta_m^*| -
(2 \log_q 2d) \sum_{n=1}^{m - \log_\theta m}
|\Delta_n^*| \\
& \geq & |\Delta_m^*| \left( 1 - \frac{(2
\log_q 2d)
\theta^{m+1-\log_\theta m}}{|\Delta_m^* |} \right) \\
& \geq & |\Delta_m^*| \left( 1 - \frac{ (2
\log_q 2d) \theta}{(\theta-1) m} \right)
\end{eqnarray*}
and therefore
\begin{eqnarray}
|\Delta_m| & \geq & |\overline{\Delta}_m| -
\left(\left(1+\left\lceil \log_q (2d)\right) \right\rceil \left(
\frac{|\overline{\Delta}_m|}{\lceil \sqrt{m} ~\rceil} +1
\right) \right)  \nonumber\\
& \geq &|\Delta_m^*| \left(1 - \frac{ (2
\log_q 2d) \theta}{(\theta-1) m} \right) \left( 1-\frac{2\left( 1
+ \left\lceil \log_q (2 d) \right\rceil \right)}{\lceil \sqrt{m}
~\rceil}\right) \label{Deltam}.
\end{eqnarray}
Also, it is clear that
\begin{equation} \label{Deltal}
\theta^{m} (\theta-1) \leq |\Delta_m^*| \leq
\theta^{m} (\theta-1) + 1,
\end{equation}
and
\begin{equation} \label{dh}
\sum_{h \geq 0} \left|\Delta_m^{(h)}\right|^3 \leq \left( \max_{h
\geq 0} \left| \Delta_m^{(h)} \right| \right)^2 \sum_{h \geq 0}
\left|\Delta_m^{(h)} \right| \leq m |\Delta_m|.
\end{equation}
By construction
$$
\sum_{h \geq 0} \sum_{k \in \Delta_m^{(h)}} p (\nu_k x) = \sum_{k
\in \Delta_m} p(\nu _k x),\qquad m \geq 1,
$$
and since we have assumed $\|p\|=1$, we also have
\begin{equation} \label{aj}
|a_j| \leq \sqrt{2}, \qquad 1 \leq j \leq d,
\end{equation}
and
\begin{equation} \label{sum}
\sum_{k \in \Delta_m} \sum_{j=1}^d \frac{a_j^2}{2} = |\Delta_m|,
\qquad m \geq 1.
\end{equation}
Finally, we have for $m \geq 1$
\begin{equation} \label{beta}
|\beta_m| \leq \sum_{h \geq 0} ~\sum_{k_1,k_2 \in \Delta_m^{(h)}}
\sum_{j_1,j_2=1}^d 1 \leq d^2 \left(\max_{h \geq 0}
|\Delta_m^{(h)}|\right) \sum_{h \geq 0} \left|\Delta_m^{(h)}
\right| \leq d^2 \sqrt{m} |\Delta_m|.
\end{equation}

\begin{lemma} \label{2dim}
Let $m,n \geq 1$, and assume that
\begin{equation} \label{mn}
m \leq n - \lceil \log_\theta n \rceil.
\end{equation}
Then for sufficiently large $m,n$ we have
$$
\left\|\varphi_{m,n} (s,t) - e^{-(s^2+t^2)/2} \right\| \leq
\frac{1}{m^4+n^4}
$$
provided
$$
|s|\le m^{1/8} \quad  |t| \le n^{1/8}.
$$
\end{lemma}

\noindent
\emph{Proof:~} Using
\begin{equation} \label{eix}
e^{i x} = (1 + i x) e^{-x^2/2 + w(x)}, \qquad |w(x)| \leq |x|^3,
\end{equation}
we have
\begin{eqnarray*}
e^{i s \eta_m} & = & \prod_{h \geq 0} \exp \left( \frac{i s \sum_{k
\in \Delta_m^{(h)}} \sum_{j=1}^d a_j \cos (2 \pi j \nu_k
x)}{\sqrt{|\Delta_m|}} \right)
\\
& = & \alpha_m(s) ~\exp \left( \sum_{h \geq 0} \frac{-s^2
\left(\sum_{k \in \Delta_m^{(h)}} \sum_{j=1}^d a_j \cos(2 \pi j
\nu_k x)\right)^2 }{2 |\Delta_m|}\right) \times \\
& & \times \exp \left( \sum_{h \geq 0} ~w \left( \frac{s \sum_{k \in
\Delta_m^{(h)}} \sum_{j=1}^d a_j \cos(2 \pi j \nu_k
x)}{\sqrt{|\Delta_m|}} \right) \right).
\end{eqnarray*}
Since
\begin{eqnarray*}
& & \sum_{h \geq 0} \left( \sum_{k \in \Delta_m^{(h)}}
\sum_{j=1}^d a_j \cos(2 \pi j \nu_k x) \right)^2 \\
& = & \sum_{h \geq 0} \sum_{k_1,k_2 \in\Delta_m^{(h)}}
\sum_{j_1,j_2=1}^d \frac{a_{j_1} a_{j_2}}{2} \left( \cos (2 \pi
(j_1 \nu_{k_1} + j_2 \nu_{k_2}) x )
+ \cos (2 \pi (j_1 \nu_{k_1} - j_2 \nu_{k_2}) x ) \right) \\
& = & \sum_{h \geq 0}~\sum_{k_1,k_2 \in \Delta_m}
\sum_{j_1,j_2=1}^d \frac{a_{j_1}
a_{j_2}}{2} \cos (2 \pi (j_1 \nu_{k_1} + j_2 \nu_{k_2}) x )\\
& & |\Delta_m| + \sum_{h \geq 0} ~\underbrace{\sum_{k_1,k_2 \in
\Delta_m} \sum_{j_1,j_2=1}^d}_{(k_1,j_1) \neq (k_2,j_2)}
\frac{a_{j_1} a_{j_2}}{2} \cos (2 \pi (j_1 \nu_{k_1} - j_2 \nu_{k_2}) x ) \\
& = & |\Delta_m| + \beta_m,
\end{eqnarray*}
where we used (\ref{sum}), we can write
\begin{eqnarray}
e^{i s \eta_m} & = & \alpha_m(s) ~\exp \left( -\frac{s^2}{2} \left(
1+ \frac{\beta_m}{|\Delta_m|}\right)+ w_m (s) \right),
\label{agr3}
\end{eqnarray}
where
\begin{eqnarray}
w_m(s) & = & \sum_{h \geq 0} w \left( \frac{s \sum_{k \in
\Delta_m^{(h)}} \sum_{j =1}^d a_j \cos(2 \pi j \nu_k
x)}{\sqrt{|\Delta_m|}} \right), \nonumber
\end{eqnarray}
and by (\ref{dh}), (\ref{aj}), (\ref{eix}), 
\begin{eqnarray}
\left| w_m(s) \right| & \leq & \left| \sum_{h \geq 0} \frac{s^3
d^3 2^{3/2} \left|\Delta_m^{(h)}\right|^3}{|\Delta_m|^{3/2}}
\right|
\nonumber\\
& \leq & 3 \left|s\right|^3 d^3 m |\Delta_m|^{-1/2}. \label{west}
\end{eqnarray}
Note further that
\begin{eqnarray}
|\alpha_m(s)| & \leq & \prod_{h \geq 0} \left(1+\frac{s^2
\left(\sum_{k \in \Delta_m^{(h)}} \sum_{j=1}^d a_j \cos(2 \pi j
\nu_k x)\right)^2}{|\Delta_m|} \right)^{1/2} \nonumber\\
& \leq & \exp \left( \sum_{h \geq 0} \frac{s^2 \left(\sum_{k \in
\Delta_m^{(h)}} \sum_{j=1}^d a_j \cos(2 \pi j \nu_k x)\right)^2}{2
|\Delta_m|}
\right)\nonumber\\
& = & \exp \left( \frac{s^2}{2} \left(1+ \frac{\beta_m}{|\Delta_m|}
\right) \right). \label{agr}
\end{eqnarray}
Finally we observe, that for $m \leq n - \lceil \log_\theta n
\rceil$
\begin{eqnarray}
& & \E \left( \alpha_m (s) \alpha_n(t) \right) \nonumber\\
& = & \int_0^1 \prod_{h \geq 0} \left(1+\frac{i s \sum_{k \in
\Delta_m^{(h)}} \sum_{j =1}^d a_j \cos(2 \pi j \nu_k
x)}{\sqrt{|\Delta_m|}} \right) \times \nonumber\\
& & \times \prod_{h \geq 0} \left(1+\frac{i t \sum_{k \in
\Delta_n^{(h)}} \sum_{j=1}^d a_j \cos(2 \pi j \nu_k
x)}{\sqrt{|\Delta_n|}}
\right) dx \nonumber\\
& = & 0, \label{agr2}
\end{eqnarray}
since by the construction of the sets $\Delta_m$, $\Delta_n$ and
$\Delta_m^{(h)}$, $\Delta_n^{(h)}$ we have for any $1 \leq j_1,
j_2 \leq d$
\begin{eqnarray*}
\frac{j_1 \nu_{k_1}}{j_2 \nu_{k_2}} \not\in \left[ \frac{1}{2},2
\right] & & \textrm{if} \quad k_1 \in \Delta_m^{(h_1)},~k_2 \in
\Delta_m^{(h_2)} \quad \textrm{for some} \quad h_1 \neq h_2 \\
& & \textrm{or} \quad k_1 \in \Delta_n^{(h_1)},~k_2 \in
\Delta_n^{(h_2)} \quad \textrm{for some} \quad h_1 \neq h_2,
\end{eqnarray*}
and also
$$
\frac{j_1 \nu_{k_1}}{j_2 \nu_{k_2}} \not\in \left[ \frac{1}{2},2
\right] \quad \textrm{if} \quad k_1 \in \Delta_m^{(h_1)},~k_2 \in
\Delta_n^{(h_2)} \quad \textrm{for some} \quad h_1,h_2 \geq 0.
$$
In fact, if e.g. $k_1 \in \Delta_m^{(h_1)}$, $k_2 \in
\Delta_m^{(h_2)}$ for $h_1 < h_2$, then necessarily
\begin{eqnarray*}
\frac{j_1 n_{k_1}}{j_2 n_{k_2}} & \leq & \frac{d
n_{k_1}}{n_{k_1+\lceil \log_q (2d) \rceil}} < d q^{\log_q (2d)}
\leq 1/2,
\end{eqnarray*}
and for $k_1 \in \Delta_m$, $k_2 \in \Delta_n$ we have
$$
\frac{j_1 n_{k_1}}{j_2 n_{k_2}} \not\in \left[ \frac{1}{2}, 2
\right] \qquad \textrm{since} \qquad \frac{n_{k_1}}{n_{k_2}}
\not\in \left[ \frac{1}{2d}, 2d \right]
$$
by (\ref{constro}) and (\ref{mn}). Thus by (\ref{agr3}), (\ref{agr}), (\ref{agr2})
\begin{eqnarray}
& & \left| \varphi_{m,n} (s,t) - e^{-s^2/2 - t^2/2} \right| \nonumber\\
& = & \left| \E \left( \alpha_m(s) \alpha_n(t)  \exp \left(
-\frac{s^2}{2} \left(1+ \frac{\beta_m}{|\Delta_m|}\right) + w_m
(s) \right) \right.\right. \times \nonumber\\
& & \left.\left. \times \exp \left( -\frac{t^2}{2} \left(1+
\frac{\beta_n}{|\Delta_n|}\right) + w_n (t) \right) \right) -
e^{-s^2/2 -
t^2/2} \right| \nonumber\\
& = & \left| \E \left( \alpha_m(s) \alpha_n(t)  \left( \exp \left(
-\frac{s^2}{2} \left(1 + \frac{\beta_m}{|\Delta_m|}\right) + w_m
(s) \right) \right.\right.\right. \times \nonumber\\
& & \left.\left.\left. \times \exp \left( -\frac{t^2}{2}\left(1 +
\frac{\beta_n}{|\Delta_n|}\right) + w_n (t) \right) - e^{-s^2/2 -
t^2/2} \right) \right) \right|
\nonumber\\
& \leq & \E \left( \left| \alpha_m(s) \alpha_n(t) \right| \left| \exp
\left( -\frac{s^2}{2} \left(1+ \frac{\beta_m}{|\Delta_m|}\right)+
w_m (s) \right) \right.\right. \times \nonumber\\
& & \left.\left. \times \exp \left( -\frac{t^2}{2}\left(1 +
\frac{\beta_n}{|\Delta_n|}\right) + w_n (t)
\right) - e^{-s^2/2 - t^2/2} \right| \right) \nonumber\\
& \leq & \E \left| e^{w_m (s)+w_n(t)} - \exp \left(\frac{s^2
\beta_m}{2 |\Delta_m|} + \frac{t^2
\beta_n}{2 |\Delta_n|} \right) \right| \nonumber\\
& \leq & \E \left| e^{w_m (s)+w_n(t)} - 1 \right| + \E \left| e
\left(\frac{s^2 \beta_m}{2 |\Delta_m|} + \frac{t^2 \beta_n}{2
|\Delta_n|} \right) - 1 \right|. \label{equ1}
\end{eqnarray}
By (\ref{west}) we have
\begin{eqnarray} \label{equ2}
\E \left| e^{w_m (s)+w_n(t)} - 1 \right| \leq e \left( 3 |s|^3 d^3
m |\Delta_m|^{-1/2} + 3 |t|^3 d^3 n |\Delta_n|^{-1/2} \right) - 1.
\end{eqnarray}
The function $\beta_m$ is a sum of at most $2 \sqrt{m} |\Delta_m|$
trigonometric functions. The coefficients of these functions are
bounded by some constant $C^*$ by the Diophantine condition in Theorem \ref{thlils}. Using (\ref{beta}), this implies
\begin{eqnarray}
\| \beta_m\|^2 & \leq & 2 C^* \sqrt{m} |\Delta_m|, \nonumber\\
\p \left( |\beta_m| > |\Delta_m|^{2/3} \right) & \leq & \frac{2
C^* \sqrt{m}}{|\Delta_m|^{-1/3}},\nonumber\\
\E ~\exp\left(\frac{s^2 \beta_m}{|\Delta_m|}\right) & \leq & \exp
\left(s^2 |\Delta_m|^{-1/3} \right) + \exp \left( s^2 d^2 \sqrt{m}
\right) \frac{2 C^* m}{|\Delta_m|^{-1/3}} \label{equ3}
\end{eqnarray}
and therefore
\begin{eqnarray}
& & \E \left| \exp \left(\frac{s^2 \beta_m}{2 |\Delta_m|} + \frac{t^2
\beta_n}{2 |\Delta_n|} \right) - 1 \right|
\nonumber\\
& \leq & \left( \E \exp \left(\frac{s^2 \beta_m}{|\Delta_m|}\right)
\right)^{1/2} \left( \E \exp \left(\frac{t^2
\beta_n}{|\Delta_n|}\right) \right)^{1/2} - 1 \label{equ4}\\
& \leq & \left( \left( \exp \left(s^2 |\Delta_m|^{-1/3} \right) + \exp
\left( s^2 d^2 \sqrt{m} \right) \frac{2 C^* m}{|\Delta_m|^{-1/3}}
\right)
\right. \times \nonumber\\
& & \left. \times ~ \left( \exp \left(t^2 |\Delta_n|^{-1/3} \right) +
\exp \left(t^2 d^2 \sqrt{n} \right) \frac{2 C^* n}{|\Delta_n|^{-1/3}}
\right) \right)^{1/2} - 1 \label{18}
\end{eqnarray}
Now (\ref{equ1}), (\ref{equ2}), (\ref{18}) and some elementary
calculations show that for sufficiently large $m,n$, under the
additional condition
$$
|s| \leq m^{1/8}, \quad |t| \leq n^{1/8}
$$
we have
$$
\left|\varphi_{m,n} (s,t) - e^{-s^2/2 - t^2/2} \right| \leq
\frac{1}{m^4+n^4},
$$
which proves the lemma. \quad $\square$

\begin{lemma} \label{1dim}
For sufficiently large $m$ we have
$$
\left|\E e^{i s \eta_m} - e^{-s^2/2} \right| \leq m^{-4},
$$
for all $s \in [-m^{1/8},m^{1/8}]$.
\end{lemma}

\noindent
\emph{Proof:~}This lemma is an one-dimensional version of Lemma
\ref{2dim} and can be shown in exactly the same way. 

\begin{lemma} \label{lemmaB}
Let $B$ be a finite set of positive integers. Then if $|B|$ is
sufficiently large, we can divide $B$ into two disjoint sets
$B_1,B_2$, such that
$$
|B_2| \leq C_1 |B|/\sqrt{\log|B|}
$$
for some constant $C_1$, and
$$
\left| \E \exp \left( i s |B_1|^{-1/2} \sum_{k \in B_1} p(n_k x) \right) - e^{-s^2/2} \right| \leq (\log |B_1|)^{-4},
$$
for $|s| \le  (\log |B|)^{1/8}$.
\end{lemma}

\emph{Proof:~}This lemma can be shown in the same way as the previous
two lemmas (or exactly in the same way as \cite[Lemma 4.3]{abt}).

\begin{lemma}[{\cite[Lemma 2]{ta}}] \label{lemmatak}
Let $B$ be a finite set of positive integers. Then for any
$\lambda>0$ satisfying
$$
4 \lambda |B|^{1/3} < 1
$$
we have
$$
\int_0^1 \exp \left( \sum_{k \in B} p(n_{\sigma(k)} x) \right) dx
\leq C_2  e^{C_3 \lambda^2 |B|},
$$
where $C_2, C_3$ are positive constants.
\end{lemma}

Next we prove some Berry-Esseen type lemmas needed for our proof. We redefine the random variables $\eta_1, \eta_2, \dots$ on a
larger probability space $(\Omega, \mathcal{A},\hat{P})$ (we write
$\hat{\eta}_1, \hat{\eta}_2,\dots$ for the redefined r.v.'s), such
that their finite dimensional distributions remain unchanged, and
such that on the new probability space there exists a sequence
$\hat{h}_1, \hat{h}_2, \dots$ of i.i.d. random variables
satisfying
\begin{itemize}
\item $\hat{h}_m \sim
\mathcal{N}(0,\tau_m),\quad \textrm{where} \quad \tau_m
=\dfrac{1}{\sqrt{8} \log \log \log \theta^m}, \quad m \geq 1$
\item $\hat{h}_m \quad \textrm{and} \quad \eta_m \quad \textrm{are
independent}, \quad m \geq 1,$
\item the two-dimensional random variables $(\hat{h}_m,\hat{h}_n)$ and
$(\hat{\eta}_m,\hat{\eta}_n)$ are independent, $m \neq n, ~m,n
\geq 1$
\end{itemize}

\begin{lemma} \label{2dimbc}
Define
$$
z_m = \sqrt{(2-\ve)(\theta/(\theta-1)) \log \log \theta^m}
$$
and
\begin{eqnarray*}
A_m & = & \left\{ \omega \in \Omega:~\hat{\eta}_m(\omega) +
\hat{h}_m (\omega)
> z_m\right\}, \quad
m \geq 1.
\end{eqnarray*}
Then
\begin{eqnarray*}
& & \left| \hat{P} (A_m A_n) - R \left( (1+\tau_m)^{-1} z_m \right)
R\left( (1 + \tau_n)^{-1} z_n \right) \right| \\
& \leq & (\log m)^2
(\log n)^2 \left( m^{-4} + n^{-4} \right)
\end{eqnarray*}
for sufficiently large $m,n$, provided $m \leq n - \lceil \log n
\rceil$. Here
$$
R(u)= 1 - (2 \pi)^{-1/2} \int_{-u}^u e^{-s^2/2} ds, \quad
u \geq 0.
$$
\end{lemma}

\emph{Proof:~} We define two measures $P_1, P_2$ on $\R^2$: $P_1$
is the measure induced by $(\hat{\eta}_m,\hat{\eta}_n)$, and $P_2$
is a two-dimensional standard normal distribution.
We apply Lemma \ref{lemmacf} with $x=z_1,y=z_2$,
$\sigma_1=\tau_m,\sigma_2=\tau_n$ and
$$
\begin{array}{ll}
T_1 = 8 \sqrt{\log \log \theta^m} \log \log \log \theta^m \qquad &
T_2 = 8 \sqrt{\log \log \theta^n} \log \log \log \theta^n.
\end{array}
$$
Then we get, using the notations from Lemma \ref{lemmacf},
\begin{eqnarray*}
& & \left| P_1^* ([-x,x] \times [-y,y]) - P_2^* \left([-x,x] \times [-y,y]\right) \right|\\
& \leq & +x y ~4 T_1 T_2 \frac{1}{m^4+n^4} \\
& & + xy  \left(  \tau_m^{-1} \tau_n~ e\left(- T_1^2
\tau_1^2/2 \right) + \tau_m \tau_n^{-1} ~e \left(-
T_2^2 \tau_2^2/2 \right) \right) \\
& \leq & (\log m)^2 (\log n)^2 \left(\frac{1}{m^4} + \frac{1}{n^4}\right)
\end{eqnarray*}
for sufficiently large $m,n$ (we emphasize that $T_1 \leq
m^{1/8},~T_2 \leq n^{1/8}$ for sufficiently large $m,n$, and
therefore we can use Lemma \ref{2dim}). Since by construction
$(\hat{\eta}_m,\hat{\eta}_n)$ and $(\hat{h}_m,\hat{h}_n)$ are
independent,
$$
\hat{P} (A_m A_n) = 1 - \left( P_1 \star H \right) ([-z_m,z_m]
\times [-z_n,z_n]) = 1 - P_1^* ([-z_m,z_m] \times [-z_n,z_n]),
$$
and since the random variables $\hat{h}_m$ have distribution
$\mathcal{N}(0,\tau_m)$, we have
$$
1 - P_2^* \left([-z_m,z_m] \times [-z_n,z_n]\right) = R\left(
(1+\tau_m)^{-1} z_m) \right) R\left( (1 + \tau_n)^{-1} z_m\right).
$$
Summarizing our estimates, we have
\begin{eqnarray*} \left| \hat{P} (A_m A_n) - R \left( (1+\tau_m)^{-1}
z_m \right) R\left( (1 + \tau_n)^{-1} z_n \right)\right| & \leq &
(\log m)^2 (\log n)^2 \left( \frac{1}{m^4} + \frac{1}{n^4} \right)
\end{eqnarray*}
for sufficiently large $m,n$. \qquad $\square$\\

\begin{lemma} \label{1dimbc}
For sufficiently large $m$
$$
\left| \hat{P} (A_m) - R \left( (1 + \tau_m)^{-1} z_m \right)
\right| \leq \frac{(\log m)^2}{m^4}.
$$
\end{lemma}

This can be shown like Lemma \ref{2dimbc}, using Lemma
\ref{lemmacf1dim} instead of Lemma \ref{lemmacf}.

\begin{lemma} \label{lemmapa}
Let
$$
\overline{A}_m = \left\{ x \in (0,1):~\sum_{k=1}^{\theta^m}
p(\nu_k x)
> \sqrt{(2+\ve) \log \log \theta^n} + 3 \frac{\sqrt{\log \log
\theta^m}}{\log \log \log \theta^m}\right\}, \quad m \geq 1.
$$
Then for sufficiently large $m$
$$
\p(\overline{A}_m) \leq R\left( \sqrt{(2+\ve) \log \log
\theta^m}\right) + 2\frac{(\log m)^2}{m^4}.
$$
\end{lemma}
\emph{Proof:~}This is a consequence of Lemma
\ref{lemmacf1dimstern} and Lemma \ref{lemmaB}. In fact, let
$$
B = \left\{ 1 \leq k \leq \theta^m \right\}.
$$
Then by Lemma \ref{lemmaB} there exist sets $B_1,B_2$ such that
such that
$$
|B_2| \leq C_1 |B|/\sqrt{\log|B|}
$$
and
$$
\left| \E \exp \left(\frac{i s \sum_{k \in B_1} p(n_k x)}{|B_1|^{1/2}}
\right) - e^{-s^2/2} \right| \leq \frac{1}{(\log |B_1|)^4},
$$
for $|s| \le (\log |B|)^{1/8}$. We apply Lemma
\ref{lemmacf1dimstern} with
\begin{eqnarray*}
T & = & 8 \sqrt{\log \log \theta^m} \log \log \log \theta^m \\
S & = & \sqrt{\log \log \theta^m} (\log \log \log \theta^m)^{-1}
\\
\sigma & = & \tau_m \\
x & = & \sqrt{(2+\ve) \log \log \theta^m} + \sqrt{\log \log
\theta^m} \, (\log \log \log \theta^m)^{-1}
\end{eqnarray*}
and get
\begin{eqnarray*}
& & \p \left\{ x \in (0,1):~\sum_{k=1}^{\theta^m} p(\nu_k x)
> \sqrt{(2+\ve) \log \log \theta^m} + 2 \frac{\sqrt{\log \log
\theta^m}}{\log \log \log
\theta^m}\right\} \\
& \leq & 1 - \left( \frac{1}{\sqrt{2 \pi}}\int_{-\sqrt{(2+\ve)
\log \log \theta^m}}^{\sqrt{(2+\ve) \log \log \theta^m}} e^{s^2/2}
~ds - 2xT \frac{1}{{(\log |B_1|)^4}} \right. \\
& & \qquad \left. - 4 x \tau_m^{-1} \exp \left( -
T^2 \tau_m^2/2 \right) - 2 \exp \left( - S^2/(2
\tau_m^2) \right) \right)
\\
& \leq & R\left( \sqrt{(2+\ve) \log \log \theta^m}\right) +
(\log m)^2 m^{-4}
\end{eqnarray*}
for sufficiently large $m$. By Lemma \ref{lemmatak}
$$
\p\left( \left| \sum_{k \in B_2} p(n_k x) \right| > S \right) \leq
m^{-4}
$$
for sufficiently large $m$, and the proof of the lemma is
complete. \qquad $\square$\\

We are ready now to prove the upper bound in the LIL. We show
\begin{lemma} \label{lemmaup}
$$
\limsup_{N \to \infty} \frac{\left| \sum_{k=1}^N p(n_{\sigma(k)}
x) \right|}{\sqrt{2 N \log \log N}} \leq 1 \quad \textup{a.e.}
$$
\end{lemma}
\emph{Proof:~} By Lemma \ref{lemmapa} we have
$$
\sum_{m \geq 1} \p \left( \overline{A}_m \right) < + \infty,
$$
and therefore the Borel-Cantelli lemma implies
\begin{equation} \label{up}
\liminf_{m \to \infty} \frac{\left| \sum_{k=1}^{\theta^m}
p(n_{\sigma(k)} x) \right|}{\sqrt{(2+\ve) \theta^m \log \log
\theta^m}} \leq 1 \quad \textup{a.e.}
\end{equation}
It remains to fill the gaps between $\theta^m$ and $\theta^{m+1}$,
$m \geq 1$. Using Lemma \ref{lemmatak} we can show, e.g. by using
the method from \cite[Section 4]{eg}, that
$$
\limsup_{m \to \infty} ~\max_{\theta^m \leq M \leq \theta^{m+1}}
\frac{\left| \sum_{k=\theta^m}^M p(n_{\sigma(k)} x)
\right|}{\sqrt{2 (\theta^{m+1}-\theta^m) \log \log
(\theta^{m+1}-\theta^m)}}  \leq C_4 \quad \textup{a.e.},
$$
where $C_4$ may only depend on $p$ and the growth factor $q$.
Combining this with (\ref{up}) we have
\begin{eqnarray*}
& & \limsup_{N \to \infty} \frac{\left|\sum_{k=1}^N p(n_{\sigma(k)} x)
\right|}{\sqrt{2 N \log \log N}} \\
& \leq & \limsup_{m \to \infty}
\max_{\theta^m \leq M \leq \theta^{m+1}} \frac{ \left|
\sum_{k=1}^M
p(n_{\sigma(k)} x) \right|}{\sqrt{2 \theta^m \log \log \theta^m}} \\
& \leq & \limsup_{m \to \infty} \frac{\left| \sum_{k=1}^{\theta^m}
p(n_{\sigma(k)} x) \right|}{\sqrt{2 \theta^m \log \log \theta^m}}
+ \limsup_{m \to \infty} \max_{\theta^m \leq M \leq \theta^{m+1}}
\frac{\left| \sum_{k=\theta^m}^M p(n_{\sigma(k)} x)
\right|}{\sqrt{2 \theta^m
\log \log \theta^m}} \\
& \leq & (2 + \ve) + C_4 (\theta-1) \quad \textup{a.e.}
\end{eqnarray*}
Since $\ve>0$ and $\theta>1$ can be chosen arbitrarily, this
concludes the proof of Lemma~\ref{lemmaup}. \qquad $\square$\\

Next we prove the lower bound in the LIL.

\begin{lemma}
$$
\limsup_{N \to \infty} \frac{\left| \sum_{k=1}^N p(n_{\sigma(k)}
x) \right|}{\sqrt{2 N \log \log N}} \geq 1 \quad \textup{a.e.}
$$
\end{lemma}
\emph{Proof:~} By Lemma \ref{2dimbc} we have
\begin{eqnarray*}
& & \sum_{n=1}^N \sum_{m=1}^N \hat{P} (A_m A_n) \\
& \geq & - C_5 + 2 \sum_{n=1}^N ~\sum_{m=n^{2/3}}^{n - \log_\theta
n} \left( R \left( (1+\tau_m)^{-1} z_m ) \right) R\left( (1 +
\tau_n)^{-1} z_n \right) \right. \\
& & \left. + (\log m)^2 (\log n)^2 \frac{1}{m^4 n^4} \right)\\
& \geq & - C_6 + 2 \sum_{n=1}^N ~\sum_{m=1}^{n} R \left(
(1+\tau_m)^{-1} z_m ) \right) R\left( (1 +
\tau_n)^{-1} z_n \right) \\
& & - 2 \sum_{n=1}^N ~\sum_{m=1}^{n^{2/3}} R \left(
(1+\tau_m)^{-1} z_m ) \right) R\left( (1 + \tau_n)^{-1}
z_n \right)\\
& & - 2 \sum_{n=1}^N ~\sum_{m=n - \log_\theta n}^{n} R \left(
(1+\tau_m)^{-1} z_m ) \right) R\left( (1 + \tau_n)^{-1} z_n
\right)
\end{eqnarray*}
for some positive constants $C_5$ and $C_6$.\\

In the sequel we will assume that $\ve$ and $\theta$ are chosen in
such a way that there exists some $\rho>0$ such that
$$
(2-\ve) (\theta/(\theta-1))/2 < 1 - \rho.
$$
For given $\ve$ this is possible by choosing $\theta$
large.
Some calculations show that
\begin{eqnarray*}
\exp \left( -\left((1 + \tau_m)^{-1} z_m +1\right)^2/2 \right) & \leq &
\sqrt{2 \pi} ~R \left( (1+\tau_m)^{-1} z_m \right)\\
& \leq & \exp \left(
-\left((1 + \tau_m)^{-1} z_m\right)^2/2 \right),
\end{eqnarray*}
and therefore
\begin{eqnarray}
& & \frac{(m \log \theta)^{-(1 +
\tau_m)^{-1}(2-\ve)(\theta/(\theta-1))/2}}{e^{-(1+\tau_m)^{-1}
\sqrt{(2-\ve)(\theta/(\theta-1))\log
\log \theta^m}/2-1/2}} \label{upb}\\
& \leq & \sqrt{2 \pi} ~R \left( (1+\tau_m)^{-1} z_m \right) \nonumber\\
& \leq & (m \log \theta)^{-(1 +
\tau_m)^{-1}(2-\ve)(\theta/(\theta-1))/2},\nonumber
\end{eqnarray}
which implies
\begin{eqnarray*}
& & \sum_{n=1}^N ~\sum_{m=1}^{n^{2/3}} R \left( (1+\tau_m)^{-1}
z_m ) \right) R\left( (1 + \tau_n)^{-1} z_n \right) \\
& & + \sum_{n=1}^N ~\sum_{m=n - \log_\theta n}^{n} R \left(
(1+\tau_m)^{-1} z_m ) \right) R\left( (1 + \tau_n)^{-1} z_n
\right) \\
& = & o \left(\sum_{n=1}^N ~\sum_{m=1}^{n} R \left(
(1+\tau_m)^{-1} z_m ) \right) R\left( (1 + \tau_n)^{-1} z_n
\right) \right) \quad \textrm{as} \quad m \to \infty.
\end{eqnarray*}
Thus
$$
\liminf_{N \to \infty} \frac{\sum_{n=1}^N \sum_{m=1}^N \hat{P}
(A_m A_n)}{\left( \sum_{m=1}^N \hat{P} (A_m) \right)^2} =1.
$$
Then, by (\ref{upb}), for sufficiently large $m$
$$
\hat{P} (A_m) \geq m^{-1+\rho/2}.
$$
Therefore
$$
\sum_{m=1}^\infty \hat{P}(A_m) = + \infty,
$$
and by Lemma \ref{rp} there occur infinitely events $A_m$ with
probability 1, which implies
$$
\limsup_{m \to \infty} \frac{\left|\hat{\eta}_m +
\hat{h}_m\right|}{\sqrt{(2-\ve)(\theta/(\theta-1)) \log \log
\theta^m}} \geq 1 \quad \textup{a.s.}
$$
Using the classical LIL for i.i.d. random variables we easily get
$$
\limsup_{m \to \infty} \frac{\left|
\hat{h}_m\right|}{\sqrt{(2-\ve)(\theta/(\theta-1)) \log \log
\theta^m}} = 0 \quad \textup{a.s.},
$$
(recall that $\tau_m \to 0$) and therefore
$$
\limsup_{m \to \infty}
\frac{\left|\hat{\eta}_m\right|}{\sqrt{(2-\ve)(\theta/(\theta-1))
\log \log \theta^m}} \geq 1 \quad \textup{a.s.}
$$
This implies the similar result for the original random variables $\eta_1, \eta_2,\dots$, i.e.
$$
\limsup_{m \to \infty}
\frac{\left|\eta_m\right|}{\sqrt{(2-\ve)(\theta/(\theta-1)) \log
\log \theta^m}} \geq 1 \quad \textup{a.e.}
$$
or
$$
\limsup_{m \to \infty} \frac{\left|\sum_{k \in \Delta_m} p(\nu_k
x)\right|}{\sqrt{(2-\ve)(\theta/(\theta-1)) |\Delta_m| \log \log
\theta^m}} \geq 1 \quad \textup{a.e.}
$$
Using Lemma \ref{lemmatak}, it is not difficult to show
$$
\limsup_{m \to \infty} \frac{\left| \sum_{k \in
\overline{\overline{\Delta}}_m \backslash \Delta_m} p(\nu_k x)
\right|}{\sqrt{(2-\ve)(\theta/(\theta-1)) |\Delta_m| \log \log
\theta^m}} = 0 \quad \textup{a.e.}
$$
and since by (\ref{Deltam}) and (\ref{Deltal})
$$
\frac{|\Delta_m|}{\theta^m (\theta-1)} \to 1
$$
this implies
\begin{eqnarray*}
\limsup_{m \to \infty} \frac{\sum_{k=\theta^m}^{\theta^{m+1}}
p(\nu_k x)}{\sqrt{(2-\ve)\theta^{m+1} \log \log \theta^{m+1}}}
\geq 1 \quad \textup{a.e.}
\end{eqnarray*}
By the results from the previous section,
$$
\limsup_{m \to \infty} \frac{\left| \sum_{k=1}^{\theta^m} p(\nu_k
x) \right|}{\sqrt{2 \theta^m \log \log \theta^m}} \leq 1 \quad
\textup{a.e.},
$$
and therefore
$$
\limsup_{m \to \infty} \frac{\sum_{k=1}^{\theta^{m+1}} p(\nu_k
x)}{\sqrt{2\theta^{m+1} \log \log \theta^{m+1}}} \geq
\frac{\sqrt{2-\ve}}{\sqrt{2}} - \frac{1}{\sqrt{\theta}} \quad
\textup{a.e.}
$$
Choosing $\ve>0$ small and $\theta>1$ large this proves
Lemma \ref{lemmap}, and therefore the proof of Theorem \ref{thlils} is complete. \qquad $\square$\\

To conclude this section,  we justify the remark made after the
statement of Theorem \ref{thlils}. Assume there exist integers $a\neq
0,b\neq 0,c$, such that the Diophantine equation
\begin{equation} \label{diophequ2}
a n_k - b n_l = c
\end{equation}
has infinitely many solutions $(k,l), k \neq l$ (by an easy
observation we can assume $a>0,b>0$). We will construct a
trigonometric polynomial $p(x)$ and a permutation $\sigma: \Z^+
\to \Z^+$ such that
\begin{equation} \label{pwert}
\limsup_{N \to \infty} \frac{\sum_{k=1}^N p(n_{\sigma(k)}
x)}{\sqrt{2 N \log \log N}} \neq \|p\| \quad \textup{a.e.}
\end{equation}
We define
$$
p(x) = \cos(2 \pi a x) + \cos(2 \pi b x).
$$
Let
$$
(k_1,l_1), (k_2,l_2), \dots
$$
denote a sequence of solutions of (\ref{diophequ2}), chosen in
such a way that 
\begin{itemize}
\item $k_j > k_i,\quad l_j > l_i \qquad \textrm{for} \ j >i$
\item $k_{i+1}/k_i > 2, \quad l_{i+1}/l_i > 2, \qquad i \geq 1$
\item $k_{i+1}/k_i \to \infty, \quad l_{i+1}/l_i
\to \infty$.
\end{itemize}
Clearly there exists a permutation $\sigma:~\Z^+ \to
\Z^+$ such that
\begin{eqnarray} \label{sigm}
\limsup_{N \to \infty} \frac{ \left| \sum_{k = 1}^N
p(n_{\sigma(k)} x) \right|}{\sqrt{2 N \log \log N}} = \limsup_{N
\to \infty} \frac{ \left| \sum_{i=1}^{N/2} p(n_{k_i} x) +
p(n_{l_i} x) \right|}{\sqrt{2 N \log \log N}}
\end{eqnarray}
For example, we can construct $\sigma$ such that for every even $N$
\begin{eqnarray*}
\left\{ \sigma(k),~1 \leq k \leq N \right\} & = & \left\{
k_i, 1 \leq i \leq N/2 - \lfloor \log_{10} N \rfloor \right\} \\
& & \cup \left\{l_i, 1 \leq i \leq N/2 - \lfloor \log_{10} N \rfloor \right\}\\
 & &
\cup \left\{ 1 \leq k \leq M \right\},
\end{eqnarray*}
where $M$ is chosen such that the set on the right-hand side
really consists of $N$ elements. Since always $M \leq 2 \log N$, relation (\ref{sigm}) will hold for $\sigma$.
Thus it suffices to calculate
$$
\limsup_{N \to \infty} \frac{ \left| \sum_{i=1}^{N/2} p(n_{k_i} x)
+ p(n_{l_i} x) \right|}{\sqrt{2 N \log \log N}}.
$$
Using standard trigonometric identities we have
\begin{eqnarray*}
& & p(n_{k_i} x) + p(n_{l_i} x) \\
& = & \cos(2 \pi a n_{k_i} x) + \cos(2 \pi b n_{k_i} x) + \cos(2
\pi a n_{l_i} x) +
\cos(2 \pi b n_{l_i} x) \\
& = & 2 \cos (\pi (a n_{k_i} + b n_{l_i}) x)\cos(\pi (a n_{k_i} -
b n_{l_i}) x) + \cos(2 \pi b n_{k_i}
x) + \cos(2 \pi a n_{l_i} x) \\
& = & 2 \cos (\pi c x)\cos\left(\pi (a n_{k_i} + b n_{l_i})
x\right) + \cos(2 \pi b n_{k_i} x) + \cos(2 \pi a n_{l_i} x).
\end{eqnarray*}
Clearly, a sequence consisting of the elements
$$
(a n_{k_i} + b n_{l_i})/2, \quad a n_{l_i}, \quad b n_{k_i},
\qquad i \geq 1,
$$
arranged in increasing order, is a lacunary sequence for $i$
sufficiently large. Using the methods of \cite{ai} we can show
\begin{eqnarray}
& & \limsup_{N \to \infty} \frac{\sum_{k=1}^N 2 \cos (\pi c
x)\cos(\pi (a n_{k_i} + b n_{l_i}) x) + \cos(2 \pi b n_{k_i} x) +
\cos(2 \pi
a n_{l_i} x)}{\sqrt{2 N \log \log N}} \nonumber\\
& = & \sqrt{2 \cos^2 (\pi c x)+1} \quad \textup{a.e.}\nonumber
\label{unbest}
\end{eqnarray}
and therefore
\begin{eqnarray*}
\limsup_{N \to \infty} \frac{ \left| \sum_{i=1}^{N/2} p(n_{k_i} x)
+ p(n_{l_i} x) \right|}{\sqrt{2 N \log \log N}} & = &
\sqrt{\cos^2 (\pi c x)+1/2} \quad \textup{a.e.} \\
& = & \sqrt{\frac{\cos (2 \pi c x)+2}{2}} \quad \textup{a.e.},
\end{eqnarray*}
which verifies (\ref{pwert}).

\end{document}